\documentclass{amsart}
\usepackage{amssymb}

\begin{document}

\newtheorem{thm}{Theorem}
\newtheorem{prop}[thm]{Proposition}
\newtheorem{lem}{Lemma}
\newtheorem*{main}{Theorem}
\newtheorem*{thecor}{Corollary}


\def \into {\hookrightarrow}
\def \onto {\twoheadrightarrow}
\def \rad {\operatorname{rad}}
\def \of {\lower3pt\hbox{${}^{\circ}$}}
\def  \Hom {\operatorname{Hom}}
\def \soc {\operatorname{soc}}
\def \gar {{\mathbb G}_{a(r)}}
\def \gas {{\mathbb G}_{a(s)}}
\def \gao {{\mathbb G}_{a(1)}}
\def \ga {{\mathbb G_a}}
\def \im {\operatorname{im}}

\title[Projectivity for unipotent group schemes]{Projectivity of modules for
infinitesimal unipotent  group schemes}
\author{Christopher P. Bendel}
\address{Department of Mathematics\\
University of Notre Dame \\
Notre Dame, IN 46556}
\curraddr{Department of Mathematics, Statistics, and Computer Science\\
University of Wisconsin-Stout\\
Menomonie,  WI  54751}
\email{bendelc@uwstout.edu}
\commby{Lance W. Small}
\begin{abstract}In this paper, it is shown that the projectivity of a rational module for an 
infinitesimal unipotent group scheme over an algebraically closed field of positive 
characteristic can be detected on a family of closed subgroups.
\end{abstract}
\maketitle

Let $k$ be an algebraically closed field of characteristic $p > 0$ and $G$ be an 
infinitesimal group scheme over $k$, that is an affine group scheme $G$ over $k$
whose coordinate (Hopf) algebra $k[G]$ is a finite-dimensional local $k$-algebra.
A rational $G$-module is equivalent to a $k[G]$-comodule and further equivalent
to a module for the finite-dimensional cocommutative Hopf algebra 
$k[G]^* \equiv \Hom_k(k[G],k)$.  Since $k[G]^*$ is a Frobenius algebra (cf.\ \cite{Jan}),
a rational $G$-module (even infinite-dimensional) is in fact projective if and only if it
is injective (cf.\ \cite{FW}).   Further,  for any rational $G$-module $M$ and any closed 
subgroup scheme $H \subset G$, if $M$ is projective over $G$, then it remains projective 
upon restriction to $H$ (cf.\ \cite{Jan}).  We consider the question of whether there is a ``nice''
collection of closed subgroups of $G$ upon which projectivity (over $G$) can be detected.

For an example of what we mean by a ``nice'' collection, consider the situation of 
modules over a finite group.  Over a field of characteristic $p > 0$, a module over a 
finite group is projective if and only if it is projective upon restriction to a 
$p$-Sylow subgroup (cf.\ \cite{Rim}).   For a $p$-group (and hence for any finite group), 
L. Chouinard \cite{Ch}  showed that a module is projective if and only if it is projective upon
restriction to every elementary abelian subgroup.   If the module is assumed to be
finite-dimensional, this result follows from the theory of varieties for finite groups (cf.\ \cite{Ca}
or \cite{Ben}).   Indeed, elementary abelian subgroups play an essential role in this theory.

In work of A. Suslin, E. Friedlander, and the author \cite{SFB1}, \cite{SFB2}, a theory of 
varieties for infinitesimal group schemes was developed.  In this setting, subgroups of the 
form $\gar$ (the $r$th Frobenius kernel of the additive group scheme $\ga$) play the role
analogous to that of elementary abelian subgroups in the case of finite groups.  Not surprisingly
then, for finite-dimensional modules, one obtains the following analogue of Chouinard's Theorem.

\begin{prop}[{\cite[Proposition 7.6]{SFB2}}]
Let $k$ be an algebraically closed field of characteristic $p> 0$, $r > 0$ be an integer, $G$ be an 
infinitesimal group scheme over $k$ of height $\leq r$, and $M$ be a finite-dimensional rational
$G$-module. Then $M$ is projective as a rational $G$-module 
if and only if whenever $H \subset G$ is a subgroup scheme isomorphic to 
$\gas$ {\rm{(}}with $s \leq r${\rm{)}} the restriction of $M$ to $H$ is 
projective as a rational $H$-module.
\end{prop}

Under some stronger hypotheses,  E. Cline, B. Parshall, and L. Scott had previously shown 
that it suffices to consider a much smaller collection of subgroups (cf.\ \cite[Main
Theorem]{CPS}).  More precisely, their result applies to a group scheme $G$ 
of the form $N_{(r)}$ (the $r$th Frobenius kernel of $N$) where
$N$ is a connected (reduced) $T$-stable subgroup scheme of a connected, semisimple 
algebraic group scheme over $k$  (and $T$ is a maximal torus).  If $M$ is a
finite-dimensional $N_{(r)}T$-module (i.e.\ $M$ is an $N_{(r)}$-module which admits 
a compatible $T$-structure), then the projectivity of $M$ may be detected by  taking 
only those subgroup schemes of the form $H = U_{\alpha (r)}$ for each root subgroup 
$U_{\alpha} \subset N$.  However, this result does not hold in general for infinite-dimensional
modules (cf.\ \cite[Example (3.2)]{CPS}), whereas Chouinard's Theorem for finite groups does.

The goal of this paper is to show that Proposition 1 also holds for infinite-dimensional
modules if $G$ is assumed to be {\em unipotent}.  (An affine group scheme $G$ is said to be
unipotent if it  admits an embedding as a closed subgroup of $U_n$, the subgroup scheme in
$GL_n$ of strictly upper triangular matrices, for some postive integer $n$.)  Although unnecessary
for Chouinard's Theorem, to handle arbitrary modules, we need a slightly stronger hypothesis. 
Indeed,  consider the case that $G$ is isomorphic to a product $\gao^{\times n}$ of height
one additive group schemes.   As an algebra, $k[G]^*$ is isomorphic to the group algebra of an
elementary abelian $p$-group and the desired theorem becomes equivalent to Dade's Lemma
\cite{Dade} which says that a $k\Pi$-module for an elementary abelian
$p$-group $\Pi$ is projective if and only if it is projective upon restriction to every ``cyclic shifted
subroup''.  This result was originally proved for finite-dimensional modules and D.\ Benson, J.\
Carlson, and J.\ Rickard \cite{BCR2} observe that one must consider a larger field in order for
Dade's Lemma to hold in general.  As such, we must also consider field extensions.  Specifically, in
Section 2, we prove the following theorem using some of the ideas in \cite{SFB2} and an argument
similar to that of Chouinard
\cite{Ch} without appealing to the theory of varieties.

\begin{main}
Let $k$ be an algebraically closed field of characteristic $p > 0$, $r > 0$ be an integer,
and $G$ be an  infinitesimal unipotent group scheme over $k$ of height $\leq r$.  For any
rational $G$-module $M$,  $M$ is projective as a rational $G$-module if and only if for
every field extension $K/k$ and every {\rm{(}}closed{\rm{)}} $K$-subgroup scheme $H \subset
G\otimes_kK$ with $H \simeq \gas\otimes_kK$ {\rm{(}}with $s \leq r${\rm{)}} the restriction of
$M\otimes_kK$ to $H$ is  projective  as a rational $H$-module.
\end{main}

The essential property of a unipotent group scheme is that the trivial module $k$ is the only
simple module or equivalently that $\Hom_G(k,M) = M^G \neq 0$ for all non-zero rational
$G$-modules $M$.  (Indeed, this property is sometimes taken as the definition of a 
unipotent group scheme.) Observe that the theory of finite-dimensional algebras shows
that if $k$ is the only simple module for such an algebra, then the algebra is indecomposable
as a module over itself, and hence a module is in fact projective if and only if it is free.
Hence, over an infinitesimal unipotent group scheme, the  notions of projective,  injective,
and free are all equivalent.

In the sense of these properties, infinitesimal unipotent group schemes are analogues of finite
$p$-groups.   For finite groups, the key to reducing the projectivity of a $G$-module
to the case of a $p$-group is that the restriction map in cohomology induced by the
embedding of a $p$-Sylow subgroup is an injection.  Unfortunately, there is no analogous
result in the setting of group schemes.  

We remind the reader that the {\em restricted} representation theory of a
restricted Lie algebra $\mathfrak{g}$ over $k$ is equivalent to the representation
theory of a certain (height 1) infinitesimal group scheme (cf.\ \cite{Jan}).  As such, the 
theorem applies to $p$-{\em nilpotent} restricted Lie algebras, which may be embedded in a 
Lie algebra of strictly upper triangular matrices, and may be stated as follows.

\begin{thecor}  Let $k$ be an algebraically closed field of characteristic $p > 0$,  
$\mathfrak{g}$ be a finite-dimensional, $p$-nilpotent restricted Lie algebra over $k$, and $M$
be a
$u(\mathfrak{g})$-module, where $u(\mathfrak{g})$ denotes the restricted enveloping algebra of
$\mathfrak{g}$.  Then $M$ is projective over $u(\mathfrak{g})$ if and only if for every field
extension $K/k$, $M$ is projective upon restriction to each subalgebra
$u(\langle x \rangle) \subset u(\mathfrak{g}\otimes_kK)$ for all $x \in 
\mathfrak{g}\otimes_kK$ with
$x^{[p]} = 0$, where $\langle x \rangle \subset \mathfrak{g}\otimes_kK$ denotes the 
one-dimensional restricted Lie subalgebra of $\mathfrak{g}\otimes_kK$ spanned by $x$, and
$x^{[p]}$ denotes the image of $x$ under the restriction map on $\mathfrak{g}\otimes_kK$.
\end{thecor}

In this context, the result of Cline, Parshall, and Scott asserts that it suffices to consider
only those elements $x$ which are root vectors.  Consider the two-dimensional nilpotent Lie
subalgebra 
$\mathfrak{g} = \langle x_{\alpha}, x_{\alpha + \beta}\rangle \subset \mathfrak{sl}_3$
generated by root vectors, where $\alpha$ and $\beta$ are the simple roots.
In Example (3.2) of \cite{CPS}, they exhibit an infinite-dimensional non-projective
$\mathfrak{g}$-module $M$ which is free upon restriction to the root vectors
$x_{\alpha}$ and $x_{\alpha + \beta}$; contradicting an infinite-dimensional
version of their result.  On the other hand, the interested reader can readily check that 
this module $M$ is not free upon restriction to the element
$x = cx_{\alpha} + dx_{\alpha + \beta}$ whenever $c, d \in k$ are
both non-zero and hence does not contradict the corollary.

Before turning to the proof, we note one potential application of the theorem 
(or preferably a generalization of it to arbitrary infinitesimal group schemes).  Recently  
D. Benson, J. Carlson, and J. Rickard \cite{BCR1}, \cite{BCR2} have developed a theory of varieties
for infinitely generated  modules over finite groups.  While Chouinard's Theorem for
finite-dimensional modules  follows from the theory of varieties, it turns out that knowing
Chouinard's Theorem in advance for arbitrary modules is necessary for certain results
in the more general theory.  Hence, an attempt to generalize the theory of varieties to arbitrary
modules for infinitesimal group schemes may in part require the above theorem.

{\bf Acknowledgments.} The author gratefully acknowledges the support of Northwestern
University.  He also thanks Eric Friedlander for the useful conversations, John
Palmieri for the discussions which led to the discovery of an error in an earlier version, and
the referee for the suggestions.

\section{Cohomology Facts}

In this section, we record two results about cohomology which will be used in the
next section to prove the theorem.  First, the following well known fact relates projectivity 
to vanishing of cohomology for unipotent group schemes. Indeed, the reader will readily 
observe that this result holds for modules over an arbitrary finite-dimensional $k$-algebra 
which admits $k$ as the only simple module and has the property that the notions of injective 
and projective are equivalent.

\begin{prop} Let $k$ be a field of characteristic $p > 0$ and $G$ be an infinitesimal unipotent 
group scheme over $k$.  For any rational $G$-module $M$, $M$ is projective 
{\rm{(}}= injective = free{\rm{)}}  if and only if $H^1(G,M) = 0$.
\end{prop}

To prove Proposition 2, we need to consider a {\em minimal} injective resolution $I_*$ for
$M$. Any rational $G$-module $M$ has a unique up to isomorphism {\em injective hull}
$Q_M$ such that $\soc Q_M \simeq \soc M$ (cf.\ \cite[I 3.17]{Jan}).  We can form
as usual an injective resolution of $M$:
$$
I_0 \overset{\delta_0}{\to} I_1 \overset{\delta_1}{\to} I_2 \to \cdots
$$ 
by taking $I_0 = Q_M$ and inductively $I_n$ to be the injective hull of 
$I_{n-1}/\im\delta_{n-2}$  and $\delta_{n-1}: I_{n-1} \onto I_{n-1}/\im\delta_{n-2} \into
I_n$ (with $I_{-1} \equiv M$ and $\delta_{-1}: M \into Q_M$).  Evidently we have
$I_{n-1}/\im\delta_{n-2} \simeq \im\delta_{n-1} = \ker\delta_n$.

\begin{lem}[{cf.\ \cite[I 2.5.4]{Ben}}] Let $I_*$ be the minimal
injective resolution of $M$ constructed above. Then the
differential in the complex $\Hom_G(k,I_*)$ is trivial and
hence we have $H^i(G,M) = \Hom_G(k,I_i) = I_i^G$ for all $i \geq 0$.
\end{lem}

\begin{proof}
Consider $\delta_i: I_i \to I_{i+1}$ for any $i \geq 0$ and let 
$\alpha: k \to I_i$ be a nonzero map.
We want to show that $\delta_i\of\alpha = 0$.
Since $\alpha$ is a nonzero map and $k$ is simple, $\alpha$ is in
fact injective and so $\alpha(k)$ is a simple submodule of $I_i$.
Hence, $\alpha(k)$ is contained in the socle of $I_i$.  By construction,
$\soc I_i = \soc(\ker\delta_i) \subset \ker\delta_i$.  Hence, 
$\alpha(k) \subset \ker\delta_i$.  In other words, $\delta_i\of\alpha = 0$
as claimed.
\end{proof}

Using Lemma 1, we now prove Proposition 2.

\begin{proof}[Proof of Proposition 2] If $M$ is projective, then clearly
$H^i(G,M) = 0$ for all $i > 0$.  Conversely, suppose that $H^1(G,M) = 0$ and let
$I_*$ be a minimal injective resolution of $M$ as above.
Lemma 1 shows in particular that $H^1(G,M) = I_1^G$ and
hence $I_1^G = 0$.  Since $G$ is unipotent, the $G$-fixed points of every non-zero 
module are non-zero.  Hence, we must have $I_1 = 0$.  But, that means $M$ was injective
(equivalently projective) to begin with.
\end{proof}

In the proof of the theorem, we will also need to make use of the
structure of the cohomology algebra of $\mathbb G_{a(1)}$, which is 
simply the same as the cohomology algebra of the finite cyclic group ${\mathbb Z}/p$.

\begin{prop}[{\cite{CPSvdK}}] If $p \neq 2$, then
the cohomology algebra $H^*(\gao,k)$ is the tensor product of a 
polynomial algebra $k[x_1]$ in one generator $x_1$ of 
degree $2$ and an exterior algebra
$\Lambda(\lambda_1)$ in one generator $\lambda_1$
of degree $1$.  If $p = 2$, then 
$H^*(\gao,k) = k[\lambda_1]$ is a polynomial
algebra in one generator $\lambda_1$ of degree $1$ {\rm{(}}and in this case
we set $x_1 = \lambda_1^2${\rm{)}}.
\end{prop}

\section{Proof of Theorem}
 We now proceed to prove the theorem.  If $M$ is projective as a rational $G$-module, then
$M\otimes_kK$ is projective over $G\otimes_kK$ and, as already noted, remains so upon
restriction to $H$.  Conversely, suppose all restrictions are projective.  Let
$k'/k$ be any (even trivial) algebraically closed field extension, $G'
= G\otimes_kk'$, and $M' = M\otimes_kk'$.  Clearly for any field extension
$K/k'$ and $K$-subgroup scheme $H \subset G'\otimes_{k'}K$ with $H
\simeq \gas\otimes_kk'\otimes_{k'}K$  the
restriction of $M'\otimes_{k'}K$ to $H$ is  projective  as a rational $H$-module.

We proceed by induction on $\dim_{k'}k'[G']$ to show the slightly stronger
statment that $M'$ is projective over $G'$ for any such $k'$.   If $\dim_{k'}k'[G']
=1$, $G'$ is the trivial group scheme and there is nothing to prove.  So, assume now
that the result holds for any $k'$ and any $H \subset G'$ with $\dim_{k'}k'[H] <
\dim_{k'}k'[G']$. If $G' \simeq \gas\otimes_kk'$ for some $s \leq r$, then we are
done by assumption.

So, we may assume from now on that $G' \not\simeq \gas\otimes_kk'$.  By
Proposition 2, to show that $M'$ is projective over $G'$, it suffices to
show that
$H^1(G',M') = 0$.   Let $\phi: G' \onto \mathbb G_{a(1)}' = \gao\otimes_kk'$ be any
non-trivial homomorphism of group  schemes.  (As $G$ is unipotent, there exists at
least one such homomorphism.)   Consider the short exact sequence of group
schemes over $k'$
$$
1 \to N \to G' \to \mathbb G_{a(1)}' \to 1
$$
and the induced Hochschild-Serre spectral sequence
$$
E_2^{p,q}(M') = H^p(\mathbb G_{a(1)}',H^q(N,M')) \Rightarrow H^{p+q}(G',M').
$$
Now, $\phi$ induces a map on cohomology 
$\phi^*: H^*(\mathbb G_{a(1)}',k') \to H^*(G',k')$.  Let $x_{\phi} = \phi^*(x_1)$
where $x_1 \in H^2(\mathbb G_{a(1)}',k')$ is the canonical generator. 
Furthermore, the spectral sequence admits an action of the cohomology ring 
$H^*(\mathbb G_{a(1)}',k')$, with $H^*(\mathbb G_{a(1)}',k')$ acting on the
abutment via $\phi^*$.  

Since $\phi$ is non-trivial, $\dim_{k'}k'[N] < \dim_{k'}k'[G']$ and so by 
induction $M'$ is projective upon restriction to $N$.  Thus,
$H^i(N,M') = 0$ for all $i > 0$ and the spectral sequence collapses to
$$
E_2^{p,0} = H^p(\mathbb G_{a(1)}',(M')^{N}) \Rightarrow H^p(G',M')
$$
giving an isomorphism $H^*(\mathbb G_{a(1)}',(M')^{N}) \simeq H^*(G',M')$.
We now recall that the action of $x_1$ induces a periodicity isomorphism
$H^i(\mathbb G_{a(1)}',Q) \simeq H^{i+2}(\mathbb G_{a(1)}',Q)$ for all $i > 0$
and any rational $\gao'$-module $Q$ (cf.\ \cite[2.3]{SFB2}).  Hence, the
action of $x_{\phi} \in H^2(G',k')$ also induces a periodicity isomorphism
$H^i(G',M') \simeq H^{i+2}(G',M')$ for all $i > 0$.

We consider the following two cases.

\noindent
{\bf CASE I:} $\dim_{k'}\Hom_{Gr/k'}(G',\mathbb G_{a(1)}') = 1$

By Theorem 1.6 of \cite{SFB2} (which is an analogue of a characterization of
$p$-groups in terms of cohomology by J.-P. Serre \cite{Se}), $x_{\phi} \in
H^2(G',k')$ is nilpotent and so we conclude from the above periodicity isomorphism
that $H^1(G',M') = 0$ and hence $M'$ is projective.

\noindent
{\bf CASE II:} $\dim_{k'}\Hom_{Gr/k'}(G',\mathbb G_{a(1)}') > 1$

Let $\phi_1$ and $\phi_2$ be two linearly independent homomorphisms
from $G'$ onto $\mathbb G_{a(1)}'$.  Further, let $c_1, c_2 \in k'$.  Then,
 by Corollary 1.5 of \cite{SFB2},
$$
c_1x_{\phi_1} + c_2x_{\phi_2} = x_{c_1^{1/p}\phi_1} + 
                                x_{c_2^{1/p}\phi_2} =
 x_{c_1^{1/p}\phi_1 + c_2^{1/p}\phi_2} \in H^2(G',k').
$$
If at least one of $c_1$, $c_2$ is non-zero, the map
$c_1^{1/p}\phi_1 + c_2^{1/p}\phi_2: G' \to \mathbb G_{a(1)}'$
is non-zero and so the action of $c_1x_{\phi_1} + c_2x_{\phi_2}$
induces a periodicity isomorphism $H^i(G',M') \simeq H^{i+2}(G',M')$ for
all $i > 0$.  In particular, $c_1x_{\phi_1} + c_2x_{\phi_2}: H^1(G',M') 
\overset{\sim}{\to} H^3(G',M')$.  Hence, if $H^1(G',M')$ was a finite-dimensional
space, then one would readily conclude from an eigenvalue argument that
$H^1(G',M') = 0$.  Since $H^1(G',M')$ might be infinite-dimensional (as $M'$ might
be),  we make use of an infinite-dimensional substitute used in \cite{BCR2}.   
For the reader's convenience, we restate this result here.

\begin{lem}[{\cite[Lemma 4.1]{BCR2}}] Let $V$ and $W$ be vector spaces over an 
algebraically closed field $k$ and $K$ be a non-trivial field extension of $k$.  Suppose
that $\phi_1, \phi_2: V \to W$ are linear maps with the property that for every pair
of scalars $c_1, c_2 \in K$, not both zero, the linear map 
$c_1\phi_1 + c_2\phi_2: V\otimes_kK \to W\otimes_kK$ is an isomorphism.  Then
both $V$ and $W$ are the zero vector space.
\end{lem}

Continuing with the proof of the theorem, let $K/k'$ be any non-trivial
algebraically closed field extension.  Consider the base changes
of $G'$ and $M'$ to 
$K$: $G_K \equiv G'\otimes_{k'}K$ and $M_K \equiv M'\otimes_{k'}K$.  Then we
have
$H^*(G',k')\otimes_{k'}K \simeq H^*(G_K,K)$ and 
$H^*(G',M')\otimes_{k'}K \simeq H^*(G_K,M_K)$.  Further, since $\phi_1$ and
$\phi_2$ remain linearly independent after base change, for any $c_1, c_2 \in K$
with at least one non-zero, the map $c_1\phi_1 + c_2\phi_2: G_K \to
\gao'\otimes_{k'}K$ is non-trivial.  Since our inductive assumption applies to any
algebraically closed field extension of $k$, it also applies to $K$ and so as above we
may conclude that 
$$
c_1x_{\phi_1} + c_2x_{\phi_2} = x_{c_1^{1/p}\phi_1 + c_2^{1/p}\phi_2} :
H^1(G_K,M_K) \overset{\sim}{\to} H^3(G_K,M_K).
$$ 
Hence, applying Lemma 2 to $K/k'$, we again conclude that 
$H^1(G',M') = 0$ and so $M'$ is projective over $G'$ and the proof is complete.


\begin{thebibliography}{CPSvdK9}

\bibitem[\bf Ben]{Ben} D.J. Benson, {\em Representations and cohomology}, 
Volumes I and II, Cambridge University Press, 1991. 

\bibitem[\bf BCR1]{BCR1}
D.J. Benson, J.F. Carlson, and J. Rickard, Complexity and varieties
for infinitely generated modules, {\em Math. Proc. Camb. Phil. Soc.}
{\bf 118} (1995), 223-243.

\bibitem[\bf BCR2]{BCR2}
D.J. Benson, J.F. Carlson, and J. Rickard, Complexity and varieties
for infinitely generated modules, II, {\em Math. Proc. Camb. Phil. Soc.}
{\bf 120} (1996), 597-615.

\bibitem[\bf Ca]{Ca} J.F. Carlson, The varieties and cohomology ring of 
a module, {\em J. Algebra} {\bf 85} (1983), 104-143. 

\bibitem[\bf Ch]{Ch} L. Chouinard, Projectivity and relative projectivity
over group rings, {\em J. Pure \& Applied Algebra} {\bf 7} (1976), 287-302.

\bibitem[\bf CPS]{CPS} E. Cline, B. Parshall, and L. Scott, On injective modules for
infinitesimal algebraic groups, I, {\em J. London Math. Soc.} (2) {\bf 31} (1985), 
277-291.

\bibitem[\bf CPSvdK]{CPSvdK} E. Cline, B. Parshall, L. Scott, and 
W. van der Kallen, Rational and generic cohomology, {\em Inventiones Math.}
{\bf 39} (1977), 143-163.

\bibitem[\bf Dade]{Dade} E.C. Dade, Endo-permutation modules over $p$-groups, II,
{\em Annals of Math.} {\bf 108} (1978), 317-346.

\bibitem[\bf FW]{FW} C.G. Faith, E.A. Walker, Direct-sum representations of
injective modules, {\em J. Algebra} {\bf 5} (1967), 203-221.

\bibitem[\bf Jan]{Jan} J.C. Jantzen, {\em Representations of Algebraic Groups},
Academic Press, 1987.

\bibitem[\bf Rim]{Rim} D.S. Rim, Modules over finite groups, {\em Annals of Math.}
{\bf 69}  (1959), 700-712.

\bibitem[\bf Se]{Se} J.-P. Serre, Sur la dimension cohomologique des groupes profinis,
{\em Topology} {\bf 3} (1965), 413-420.

\bibitem[\bf SFB1]{SFB1} A. Suslin, E.M. Friedlander, C.P. Bendel, 
Infinitesimal 1-parameter subgroups and cohomology, {\em Jour. Amer.
Math. Soc.} {\bf 10} (1997), 693-728.
 
\bibitem[\bf SFB2]{SFB2} A. Suslin, E.M. Friedlander, C.P. Bendel, 
Support varieties for infinitesimal group schemes, {\em Jour. Amer.
Math. Soc.} {\bf 10} (1997), 729-759.


\end{thebibliography}
\end{document}